\documentclass[12pt]{amsart}

\usepackage{amsmath,amssymb,latexsym,amsthm}

\newcommand{\des}{{\rm des}}

\newcommand{\sumlim}{\sum\limits}
\newcommand{\prodlim}{\prod\limits}

\newtheorem{thm}{Theorem}[section]

\newtheorem{lem}[thm]{Lemma}

\newtheorem{cor}[thm]{Corollary}

\newtheorem{obs}[thm]{Observation}

\newtheorem{exa}[thm]{Example}
\newtheorem{rem}[thm]{Remark}

\theoremstyle{definition}
\newtheorem{defn}[thm]{Definition}

\newtheorem{prop}[thm]{Proposition}
\newtheorem{conj}[thm]{Conjecture}
\newtheorem{clm}[thm]{Claim}
\newcommand{\een}{\end{enumerate}}
\newcommand{\blem}{\begin{lem}}
\newcommand{\elem}{\end{lem}}
\newcommand{\bcl}{\begin{clm}}
\newcommand{\ecl}{\end{clm}}
\newcommand{\bthm}{\begin{thm}}
\newcommand{\ethm}{\end{thm}}
\newcommand{\bpr}{\begin{prop}}
\newcommand{\epr}{\end{prop}}
\newcommand{\bco}{\begin{cor}}
\newcommand{\eco}{\end{cor}}
\newcommand{\bcon}{\begin{conj}}
\newcommand{\econ}{\end{conj}}
\newcommand{\bde}{\begin{defn}}
\newcommand{\ede}{\end{defn}}
\newcommand{\bex}{\begin{exa}}
\newcommand{\eexa}{\end{exa}}
\newcommand{\bobs}{\begin{obs}}
\newcommand{\eobs}{\end{obs}}
\newcommand{\bexe}{\begin{exe}}
\newcommand{\eexe}{\end{exe}}

\newcommand{\grn}{G_{r,n}}
\newcommand{\exc}{{\rm exc}}
\newcommand{\negg}{{\rm neg}}
\newcommand{\csum}{{\rm csum}}

\begin{document}

\title[Recursions for excedances in permutations groups]{Recursions for Excedance number in some permutations groups}

\author[Bagno]{Eli Bagno}
\address{The Jerusalem College of Technology, Jerusalem, Israel}
\email{bagnoe@jct.ac.il}

\author[Garber]{David Garber}
\address{Department of Applied Mathematics, Faculty of Sciences, Holon Institute of Technology, PO Box 305,
58102 Holon, Israel} \email{garber@hit.ac.il}

\author[Mansour]{Toufik Mansour}
\address{Department of Mathematics, University of Haifa, 31905
Haifa, Israel} \email{toufik@math.haifa.ac.il}

\author[Shwartz]{Robert Shwartz}
\address {Department of Mathematics, Bar-Ilan University, 52900 Ramat-Gan, Israel}
\email{shwartr1@macs.biu.ac.il}

\begin{abstract}

The excedance number for $S_n$ is known to have an Eulerian
distribution. Nevertheless, the classical proof uses descents rather
than excedances. We present a direct recursive proof which seems to
be folklore and extend it to the colored permutation groups
$G_{r,n}$.

The generalized recursion yields some interesting connection to
Stirling numbers of the second kind. We also show some log-concavity
result concerning a variant of the excedance number.

Finally, we show that the generating function of the excedance
number defined on $G_{r,n}$ is symmetric.

\end{abstract}

\date{\today}
\maketitle

\section{Introduction}

Let $S_n$ be the symmetric group on $n$ letters. The parameter {\it
excedance}, which is defined on a permutation $\pi \in S_n$ by
$$\exc(\pi)=|\{i\in [n] \mid \pi(i)>i \}|,$$ is well-known. Another
classical parameter defined on permutations of $S_n$ is the {\it
descent number}, defined by $$\des(\pi)=|\{i\in [n-1] \mid
\pi(i)>\pi(i+1) \}|.$$ Both parameters  have the same distribution,
which can be read from the following recursion:
$$a(n,k)= (k+1) a(n-1,k) +(n-k)a(n-1,k-1);$$
$$\quad a(n,0)=1,\ a(0,k)=0$$ where $a(n,k)$ is the number of
permutations in $S_n$ with $k$ excedances or $k$ descents. The
corresponding generating function $\sumlim_{\pi \in
S_n}q^{\exc(\pi)}$ is called {\it the Eulerian polynomial}.

There is a well-known proof for this recursion by enumerating the
descents \cite{FS}, and there is a bijection from $S_n$ onto
itself, taking the descents into the excedances \cite{Sta}.

We start this paper by presenting a classical way to obtain the
above recursion using only the excedance numbers. This argument
appears also in \cite{J}.

An excedance number was defined also for the family of groups
$G_{r,n}=\mathbb{Z}_r \wr S_n$ (see \cite{BG},\cite{F},\cite{St}).

We generalize this recursion for the cases of the hyperoctahedral
group $B_n=G_{2,n}$ and the colored permutation groups $G_{r,n}$.

This generalized recursion yields several interesting results.
First, it gives some connection to the Stirling numbers of the
second kind. Moreover, one gets a log-concavity result for a variant
of the excedance number.

It is well-known that the generating function $\sumlim_{\pi \in
S_n}{{q^{\exc(\pi)}}}=\sumlim_{i=0}^d{a_iq^i}$ has some symmetric
properties. It is symmetric in the sense that $a_i=a_{d-i}$ for $i
\in \{1,\cdots,\lfloor \frac{d}{2} \rfloor\}$. We prove here the
corresponding symmetric property for $\grn$.

\medskip

The paper is organized as follows. In Section \ref{grn}, we
introduce the group $G_{r,n} = \mathbb{Z}_r \wr S_n$. In Section
\ref{stat}, we define some of its parameters.

Section \ref{proof_Sn} deals with the proof of the recursion for
$S_n$. In Sections \ref{proof_Bn} and \ref{proof_Grn} we give the
corresponding recursions for $B_n$ and $G_{r,n}$, respectively.

Section \ref{exc_A_dist} deals with the connection to Stirling
numbers, and Section \ref{concave} deals with the log-concavity
result.

In the last section, we prove the symmetry of the generating
function of the excedance on $G_{r,n}$.

\section{The group of colored permutations}\label{grn}

\bde Let $r$ and $n$ be positive integers. {\it The group of colored
permutations of $n$ digits with $r$ colors} is the wreath product:
$$\grn=\mathbb{Z}_r \wr S_n=\mathbb{Z}_r^n \rtimes S_n,$$
consisting of all the pairs $(\vec{z},\tau)$ where
$\vec{z}$ is an $n$-tuple of integers between $0$ and $r-1$ and $\tau \in
S_n$. The multiplication is defined by the following rule: for
$\vec{z}=(z_1,\dots,z_n)$ and $\vec{z}'=(z'_1,\dots,z'_n)$
$$(\vec{z},\tau) \cdot (\vec{z}',\tau')=((z_1+z'_{\tau^{-1}(1)},\dots,z_n+z'_{\tau^{-1}(n)}),\tau \circ \tau')$$ (here $+$ is
taken modulo $r$). \ede

We use some conventions along this paper. For an element
$\pi=(\vec{z},\tau) \in \grn$ with $\vec{z}=(z_1,\dots,z_n)$ we
write $z_i(\pi)=z_i$. For $\pi=(\vec{z},\tau)$, we denote
$|\pi|=(\vec{0},\tau)$ where $\vec{0} \in \mathbb{Z}_r^n$. We also
define $c_i(\pi)=z_i(\pi^{-1})$ and $\vec{c}=(c_1,\dots,c_n)$. Using
this notation, the element
$(\vec{c},\tau)=\left((0,1,2,3),\begin{pmatrix} 1 & 2 & 3 &4 \\
2 & 1 & 4 &3
\end{pmatrix}\right) \in
G_{3,4}$ will be written as $(2 \bar{1} \bar{\bar{4}}
\bar{\bar{\bar{3}}}).$

Here is another way to present $\grn$: Consider the alphabet
$\Sigma=\{1,\dots,n,\bar{1},\dots,\bar{n},\dots,
1^{[r-1]},\dots,n^{[r-1]} \}$ as the set $[n]$ colored by the colors
$0,\dots,r-1$. Then, an element of $\grn$ is a {\it colored
permutation}, i.e., a bijection $\pi: \Sigma \rightarrow \Sigma$
satisfying the following condition: if
$\pi(i^{[\alpha]})=j^{[\beta]}$ then
$\pi(i^{[\alpha+1]})=j^{[\beta+1]}$.

In particular, $G_{1,n}=C_1 \wr S_n$ is the symmetric group $S_n$,
while $G_{2,n}=C_2\wr S_n$ is the group of signed permutations
$B_n$, also known as the {\it hyperoctahedral group}, or the {\it
classical Coxeter group of type B}.

\section{Statistics on $\grn$}\label{stat}

Given any ordered alphabet $\Sigma'$, we recall the definition of
the {\it excedance set} of a permutation $\pi$ on $\Sigma'$:
$${\rm Exc}(\pi)=\{i \in \Sigma' \mid \pi(i)>i\}$$ and the {\it excedance
number} is defined to be ${\rm exc}(\pi)=|{ \rm Exc}(\pi)|$.
\smallskip

We start by defining the color order on the set
$$\Sigma=\{1,\dots,n,\bar{1},\dots,\bar{n},\dots,
1^{[r-1]},\dots,n^{[r-1]} \}.$$

\bde The {\it color order} on $\Sigma$ is defined to be
$$1^{[r-1]} <\! \cdots < n^{[r-1]} < 1^{[r-2]} < 2^{[r-2]} < \cdots < n^{[r-2]} < \cdots < 1 < \cdots\! < n.$$
\ede

\begin{exa}

Given the color order
$$\bar{\bar{1}} < \bar{\bar{2}}
<\bar{\bar{3}} < \bar{1} < \bar{2} <\bar{3} < 1 < 2 < 3,$$ we
write $\sigma=(3\bar{1}\bar{\bar{2}}) \in G_{3,3}$ in an extended
form,
$$\begin{pmatrix} \bar{\bar{1}} & \bar{\bar{2}} & \bar{\bar{3}} &
\bar{1} & \bar{2}& \bar{3} & 1 & 2 & 3\\
\bar{\bar{3}} & 1 & \bar{2} & \bar{3} & \bar{\bar{1}}  &  2 & 3 &
\bar{1} & \bar{\bar{2}}
\end{pmatrix}$$
and calculate ${\rm
Exc}(\sigma)=\{\bar{\bar{1}},\bar{\bar{2}},\bar{\bar{3}},\bar{1},\bar{3},1\}$
and ${\rm exc}(\sigma)=6$.
\end{exa}

We present now an alternative way to compute the excedance number.
Let $\sigma \in \grn$. We define:
$${\rm csum}(\sigma) = \sumlim_{i=1}^n c_i(\sigma).$$

Note that in the case $r=2$ (i.e. the group $B_n$) the alphabet
$\Sigma$ can be seen as containing the digits $\{\pm 1,\dots ,\pm
n\}$ and the parameter $\csum(\pi)$ counts the number of digits $i
\in [n]$ such that $\pi(i)<0$, so we call it $\negg(\pi)$.

Define now:

$${\rm Exc}_A(\sigma) = \{ i \in [n-1] \ | \ \sigma(i) > i \},
$$
where the comparison is with respect to the color order, and
$${\rm exc}_A(\sigma) = |{\rm Exc}_A(\sigma)|.$$

\begin{exa}
Take $\sigma=(\bar{1}\bar{\bar{3}}4\bar{2}) \in G_{3,4}$. Then
${\rm csum}(\sigma)=4$,\break ${\rm Exc_A}(\sigma)=\{3\}$ and
hence ${\rm exc}_A(\sigma)=1$.
\end{exa}

We have now (see \cite{BG}):

\begin{lem}
\begin{align*}
{\rm exc}(\sigma)&=r \cdot {\rm exc}_A(\sigma)+{\rm
csum}(\sigma).
\end{align*}
\end{lem}

\section{The recursion for $S_n$}\label{proof_Sn}

We supply a classical proof for the recursion for the Eulerian
polynomial using its interpretation as a generating function for the
excedance number for $S_n$. Denote by $a(n,k)$ the number of
permutations in $S_n$ with exactly $k$ excedances. Then we have the
following recursion:

\bpr

$$a(n,k)=(k+1)a(n-1,k)+(n-k)a(n-1,k-1)$$
$$a(n,0)=1, a(0,k)=0.$$
\epr

\begin{proof}
For each $n$, $0 \leq k \leq n-1$, denote by $S(n,k)$ the set of
permutations in $S_n$ with exactly $k$ excedances. Denote also

$$R=\{\pi \in S(n,k) \mid \pi^{-1}(n)< \pi(n)\}$$ and
$$T=\{\pi \in S(n,k) \mid \pi^{-1}(n) \geq \pi(n)\}.$$

Define $\Phi : S(n,k) \mapsto S(n-1,k) \cup S(n-1,k-1)$ as follows:
Let $\pi \in S(n,k)$. Then $\Phi(\pi)$ is the permutation of
$S_{n-1}$ obtained from $(n,\pi(n))\pi$ by ignoring the last digit.

Let $\pi \in S(n,k)=R \cup T$. If $\pi \in R$ then $\Phi(\pi) \in
S(n-1,k)$. Note that $|\Phi^{-1}(\Phi(\pi))|=k+1$. On the other
hand, if $\pi \in T$, then $\Phi(\pi) \in S(n-1,k-1)$ and
$|\Phi^{-1}(\Phi(\pi))|=n-1-(k-1)=n-k$.
\end{proof}

We give the following example for clarifying the proof.

\begin{exa}
Consider $S(5,2)$.

Let $$R \ni \pi = \begin{pmatrix} \textcircled{1} & \textcircled{2}
& 3 & 4 & 5 \\ 5 & 3 & 1 & 2 & 4
\end{pmatrix} \mapsto \begin{pmatrix} \textcircled{1} &
\textcircled{2} & 3 & 4 & 5 \\ \mathbf{4} & 3 & 1 & 2 & \mathbf{5}
\end{pmatrix} \mapsto \begin{pmatrix} \textcircled{1} &
\textcircled{2} & 3 & 4  \\ 4 & 3 & 1 & 2
\end{pmatrix},$$

so: $$\Phi^{-1}(\Phi(\pi))=\left\{\begin{pmatrix} \textcircled{1} &
\textcircled{2} & 3 & 4 & 5 \\ 5 & 3 & 1 & 2 & 4
\end{pmatrix},\begin{pmatrix} \textcircled{1} &
\textcircled{2} & 3 & 4 & 5 \\ 4 & 5 & 1 & 2 & 3
\end{pmatrix},\begin{pmatrix} \textcircled{1} &
\textcircled{2} & 3 & 4 & 5 \\ 4 & 3 & 1 & 2 & 5
\end{pmatrix}\right\}.$$

Let $$T \ni \pi=\begin{pmatrix} \textcircled{1} & 2 &
\textcircled{3} & 4 & 5
\\ 5 & 2 & 4 & 3 & 1
\end{pmatrix} \mapsto \begin{pmatrix} 1 & 2 &
\textcircled{3} & 4 & 5
\\ \mathbf{1} & 2 & 4 & 3 & \mathbf{5}
\end{pmatrix} \mapsto \begin{pmatrix} 1 & 2 &
\textcircled{3} & 4
\\ 1 & 2 & 4 & 3
\end{pmatrix}.$$

Then $$\Phi^{-1}(\Phi(\pi))=\left\{\begin{pmatrix} \textcircled{1} &
2 & \textcircled{3} & 4 & 5
\\ 5 & 2 & 4 & 3 & 1
\end{pmatrix},\begin{pmatrix} 1 & \textcircled{2} &
\textcircled{3} & 4 & 5
\\ 1 & 5& 4 & 3 & 2
\end{pmatrix},\begin{pmatrix} 1 & 2 &
\textcircled{3} & \textcircled{4} & 5
\\ 1 & 2 & 4 & 5 & 3
\end{pmatrix}\right\}.$$
\end{exa}

\section{The recursion for $B_n$}\label{proof_Bn}

In this section, we generalize the above recursion to
$B_n=\mathbb{Z}_2 \wr S_n$. We start with some notations.
$$C_i(n,k)=\{ \pi \in B_n \mid \exc_A(\pi)=k, \negg(\pi)=i \}$$
$$c_i(n,k)=\# C_i(n,k)$$
$$B_i(n,k)=\{ \pi \in B_n \mid \exc(\pi)=k, \negg(\pi)=i \}$$
$$b_i(n,k)=\# B_i(n,k)$$
$$B(n,k)=\bigcup_{i=0}^n B_i(n,k)$$
$$b(n,k)=\sumlim_{i=0}^n b_i(n,k)$$

In $B_n$, we have $\exc(\pi)=2\exc_A(\pi)+\negg(\pi)$, hence:
$$b_i(n,k)=c_i \left( n,\frac{k-i}{2} \right)$$

In the following proposition, we give a recursion for $c_i(n,k)$:
\bpr\label{Cink}
\begin{eqnarray*}
c_i(n,k)&=&(n-k)c_i(n-1,k-1)+(k+1)c_i(n-1,k)+ \\
& &+(n-k)c_{i-1}(n-1,k)+(k+1)c_{i-1}(n-1,k+1),
\end{eqnarray*}
with the following initial conditions:
$$c_i(n,0) =
\sumlim_{\begin{array}{c}(t_1,\dots,t_i) \\1 \leq  t_1 < t_2 <
\cdots <t_i \leq n
\end{array}} i! (i+1)^{n-t_i} \prodlim_{u=1}^i u^{t_u-t_{u-1}-1}$$
where $t_0=0$, and
$$c_i(0,k)=0;\ c_0(1,0)=1;$$
$$c_{-1}(n,k)=0\ \forall n \forall k$$
\epr

\begin{proof}
We start with the recursion. Define:
$$C_i^0(n,k) =\{ \pi \in B_n \mid \exc_A(\pi)=k, \negg(\pi)=i, \pi^{-1}(n)>0\}$$
$$c_i^0(n,k)=\#C_i^0(n,k)$$
$$C_i^1(n,k) =\{ \pi \in B_n \mid \exc_A(\pi)=k, \negg(\pi)=i, \pi^{-1}(n)<0\}$$
$$c_i^1(n,k)=\#C_i^1(n,k)$$

\medskip

Obviously, $C_i(n,k)=C_i^0(n,k) \cup C_i^1(n,k)$ and hence:
$$c_i(n,k)=c_i^0(n,k)+c_i^1(n,k).$$

Define $\Phi : C_i^0(n,k) \mapsto C_i(n-1,k) \cup C_i(n-1,k-1)$ as
follows: Let $\pi \in C_i^0(n,k)$. Then $\Phi(\pi)$ is the
permutation of $B_{n-1}$ obtained from $(n,\pi(n))\pi$ by ignoring
the last digit.

Now, define:
$$R_0=\{\pi \in C_i^0(n,k) \mid \pi^{-1}(n)< \pi(n)\}$$ and
$$T_0=\{\pi \in C_i^0(n,k) \mid \pi^{-1}(n) \geq \pi(n)\}.$$

Let $\pi \in C_i(n,k)=R_0 \cup T_0$. If $\pi \in R_0$, then
$\Phi(\pi) \in C_i(n-1,k)$. Note that $|\Phi^{-1}(\Phi(\pi))|=k+1$.
On the other hand, if $\pi \in T_0$, then $\Phi(\pi) \in
C_i(n-1,k-1)$ and $|\Phi^{-1}(\Phi(\pi))|=n-1-(k-1)=n-k$.

\medskip

Define $\Phi : C_i^1(n,k) \mapsto C_{i-1}(n-1,k) \cup
C_{i-1}(n-1,k+1)$ as before.

Now, define:
$$R_1=\{\pi \in C_i^1(n,k) \mid |\pi^{-1}(n)|< \pi(n)\}$$ and
$$T_1=\{\pi \in C_i^1(n,k) \mid |\pi^{-1}(n)| \geq \pi(n)\}.$$

Let $\pi \in C_i^1(n,k)=R_1 \cup T_1$. If $\pi \in R_1$, then
$\Phi(\pi) \in C_{i-1}(n-1,k+1)$. Note that
$|\Phi^{-1}(\Phi(\pi))|=k+1$. On the other hand, if $\pi \in T_1$,
then $\Phi(\pi) \in C_{i-1}(n-1,k)$ and
$|\Phi^{-1}(\Phi(\pi))|=n-1-(k-1)=n-k$.

Combining together all the parts, we get the requested recursion for
$c_i(n,k)$.

\medskip

Now we prove the initial condition:
$$c_i(n,0) =
\sumlim_{\begin{array}{c}(t_1,\dots,t_i) \\ 1\leq t_1 < t_2 < \cdots
<t_i \leq n
\end{array}} i! (i+1)^{n-t_i} \prodlim_{u=1}^i u^{t_u-t_{u-1}-1}$$
where $t_0=0$.

Let $\pi \in B_n$ and let $1 \leq t_1 <\cdots < t_i \leq n$ be such
that $\pi(t_j)<0$ ($1 \leq j \leq i$).

In order to ensure that $\exc_A(\pi)=0$, we have to demand that for
each $\ell \not\in  \{t_1,\dots, t_i \}$, $\pi(\ell) \leq \ell$. For
each $\ell < t_1$ (if there is any), we have only one possibility:
$\pi(\ell)=\ell$. For $t_1 < \ell < t_2$ (if there is any), we have
exactly two possibilities, and so on: for $t_m < \ell < t_{m+1}$ (if
there is any), we have exactly $m+1$ possibilities. Finally, for
$t_i < \ell$, we have exactly $i+1$ possibilities.

After fixing $\pi(\ell)$ for each $\ell \not\in  \{t_1,\dots, t_i
\}$, we have exactly $i!$ possibilities to locate $\pi(t_j)$, $1
\leq j \leq i$. This gives us the required initial condition.
\end{proof}

The following example should clarify the above proof of the initial
condition. Let $\pi \in B_{9}$ and assume that $t_1=3,t_2=6, t_3=8$.
Then in order to get $\exc_A(\pi)=0$, we must have
$\pi(1)=1,\pi(2)=2$. $\pi(4)$ can be $3$ or $4$. $\pi(5) \in
\{3,4,5\}$ but once $\pi(4)$ has been chosen we have only $2$
possibilities for it. $\pi(7) \in \{3,4,5,6,7\}$  which gives us $3$
possibilities and for $\pi(9)$ we have $4$ possibilities. The values
corresponding to $\{\pi(3),\pi(6),\pi(8)\}$ are already fixed so we
just have to order them.

\section{The corresponding recursion for $G_{r,n}$}\label{proof_Grn}

The recursion for $B_n$ can be generalized to
$G_{r,n}=\mathbb{Z}_r \wr S_n$ very easily. We continue with
similar notations.
$$C_i(r,n,k)=\{ \pi \in G_{r,n} \mid \exc_A(\pi)=k, \csum(\pi)=i \}$$
$$c_i(r,n,k)=\# C_i(r,n,k)$$
$$B_i(r,n,k)=\{ \pi \in G_{r,n} \mid \exc(\pi)=k, \csum(\pi)=i \}$$
$$b_i(r,n,k)=\# B_i(r,n,k)$$
$$B(r,n,k)=\bigcup_{i=0}^n B_i(r,n,k)$$
$$b(r,n,k)=\sumlim_{i=0}^n b_i(r,n,k)$$

In $G_{r,n}$, we have $\exc(\pi)=r \cdot \exc_A(\pi)+\csum(\pi)$,
hence:
$$b_i(r,n,k)=c_i \left(r, n,\frac{k-i}{r} \right)$$

In the following proposition, we give a recurrence for
$c_i(r,n,k)$:
\bpr\label{pro71}
\begin{eqnarray*}
c_i(r,n,k)&=&(n-k)c_i(r,n-1,k-1)+(k+1)c_i(r,n-1,k)+ \\
& &+\sumlim_{j=1}^{r-1}
\left((n-k)c_{i-j}(r,n-1,k)+(k+1)c_{i-j}(r,n-1,k+1)\right),
\end{eqnarray*}
with the following initial conditions:
$$c_i(r,n,0) =
\sumlim_{\begin{array}{c}(t_1,\dots,t_i) \\ 1 \leq t_1 < t_2 <
\cdots <t_i \leq n
\end{array}} i! (r-1)^i (i+1)^{n-t_i} \prodlim_{u=1}^i u^{t_u-t_{u-1}-1}$$
where $t_0=0$, and
$$c_i(r,0,k)=0;\ c_0(r,1,0)=1;$$
$$c_{-1}(r,n,k)=0\ \forall n \forall k$$
\epr

\section{The distribution of ${\rm exc}_A$}\label{exc_A_dist}

The parameter ${\rm exc}_A$ deserves a special treatment, since its
distribution involves the Stirling number of the second kind.

 Let $d(r,n,k)=\sumlim_{i=0}^{(r-1)n}c_{i}(r,n,k)$ be the number of
permutations in $ \pi \in \grn$ such that ${\rm exc}_A(\pi)=k$.

Summing up the equation of Proposition \ref{pro71} for all $i$, we
get:
$$\begin{array}{ll}
d(r,n,k)&=(n-k)d(r,n-1,k-1)+(k+1)d(r,n-1,k)\\
&+(n-k)(r-1)d(r,n-1,k)+(k+1)(r-1)d(r,n-1,k+1), \end{array}$$ which
is equivalent to:
\begin{equation}\label{eqdd1}
\begin{array}{ll}
d(r,n,k)&=(n-k)d(r,n-1,k-1)\\
&+(k+1+(r-1)(n-k))d(r,n-1,k)\\
&+(k+1)(r-1)d(r,n-1,k+1).\end{array}
\end{equation}
In order to solve this recurrence, we define the following
polynomial:
$$D_{r,n}(t)=\sum_{k=0}^nd(r,n,k)t^k.$$
Rewriting Equation \eqref{eqdd1} in terms of the polynomial
$D_{r,n}(t)$, we obtain that:
$$\begin{array}{ll}
D_{r,n}(t)&=ntD_{r,n-1}(t)-t\frac{\partial}{\partial
t}(tD_{r,n-1}(t))\\
&+(1+(r-1)n)D_{r,n-1}(t)-(r-2)t\frac{\partial}{\partial
t}D_{r,n-1}(t)\\
& +(r-1)\frac{\partial}{\partial t}D_{r,n-1}(t), \end{array}$$ which
implies that:
\begin{equation}\label{eqdd2}
D_{r,n}(t)=(rn+(n-1)(t-1))D_{r,n-1}(t)-(t-1)(t+r-1)\frac{\partial}{\partial
t}D_{r,n-1}(t).
\end{equation}
Now, in order to solve this recurrence, assume $D_{r,n}(t)$ can be
written as $D_{r,n}(t)=P_{r,n}(t)E_{r,n}(t)$, and later we give a
condition for $P_{r,n}(t)$. Therefore, Equation \eqref{eqdd2} can be
written in terms of $P_{r,n}(t)$ and $E_{r,n}(t)$ as:
\begin{eqnarray*}
E_{r,n}(t)&=&\frac{(rn+(n-1)(t-1))P_{r,n-1}(t)-(t-1)(t+r-1)\frac{\partial}{\partial
t}P_{r,n-1}(t)}{P_{r,n}(t)}E_{r,n-1}(t)\\
& &
-\frac{(t-1)(t+r-1)P_{r,n-1}(t)}{P_{r,n}(t)}\frac{\partial}{\partial
t}E_{r,n-1}(t).
\end{eqnarray*}
Let us assume that:
$$(rn+(n-1)(t-1))P_{r,n-1}(t)=(t-1)(t+r-1)\frac{\partial}{\partial
t}P_{r,n-1}(t).$$ One solution of the above differential equation is
\begin{equation}\label{eqdd3}
P_{r,n}(t)=\frac{(t-1)^{n+1}}{t+r-1}.
\end{equation}
Note that $\frac{P_{r,n-1}(t)}{P_{r,n}(t)}=\frac{1}{t-1}$.
Therefore, for all $n\geq1$,
\begin{equation}\label{eqdd4}
E_{r,n}(t)=-(t+r-1)\frac{\partial}{\partial t}E_{r,n-1}(t).
\end{equation}
Substituting $n=1$, we have:
$$E_{r,1}(t)=-(t+r-1)\frac{\partial}{\partial t}E_{r,0}(t),$$
On the other hand, note that $D_{r,1}=r$ and
$P_{r,1}=\frac{(t-1)^2}{t+r-1}$, and therefore:
$$E_{r,1}= \frac{D_{r,1}}{P_{r,1}}=\frac{r(t+r-1)}{(t-1)^2}.$$
Hence:
$$\frac{\partial}{\partial t}E_{r,0}(t)=-\frac{r}{(t-1)^2}.$$
Hence, we define $E_{r,0}(t)=\frac{r}{t-1}$.

Before stating the main result of this section, we introduce the
Stirling numbers of the second kind. Let $S_{n,j}$ be the
$(n,j)-$Stirling numbers, and recall that for all $n,j$:
$$S_{n+1,j}=jS_{n,j}+S_{n,j-1}$$

\bthm\label{thmdd5} For all $n\geq1$,
$$E_{r,n}(t)=(-1)^{n-1}r\sum_{j=1}^{n}j!S_{n,j}\frac{(t+r-1)^j}{(1-t)^{j+1}},$$
\ethm

\begin{proof}
We start by proving by induction on $n$ the following equation:
$$E_{r,n}(t)=(-1)^n\sum_{j=1}^nS_{n,j}(t+r-1)^j\frac{\partial^j}{\partial
t^j}E_{r,0}(t).$$ It is easy to verify the equation for $n=1$.
Assume its correctness for a given $n$. By the recurrence relation
\eqref{eqdd4}, we have:
$$E_{r,n+1}(t)=-(t+r-1)\frac{\partial}{\partial t}E_{r,n}(t).$$
By the induction hypothesis:
\begin{tiny}
\begin{eqnarray*}
E_{r,n+1}(t)&=&(-1)^{n+1}(t+r-1)\sum_{j=1}^nS_{n,j}\frac{\partial}{\partial
t}\left((t+r-1)^j\frac{\partial^j}{\partial t^j}E_{r,0}(t)\right)=\\
&=&(-1)^{n+1}(t+r-1)\sum_{j=1}^nS_{n,j}\left(
j(t+r-1)^{j-1}\frac{\partial^j}{\partial
t^j}E_{r,0}(t)+(t+r-1)^j\frac{\partial^{j+1}}{\partial
t^{j+1}}E_{r,0}(t) \right)=\\
&=&(-1)^{n+1} \left( \sum_{j=1}^n j S_{n,j}
(t+r-1)^j\frac{\partial^j}{\partial t^j}E_{r,0}(t) + \sum_{j=1}^n j
S_{n,j} (t+r-1)^{j+1}\frac{\partial^{j+1}}{\partial
t^{j+1}}E_{r,0}(t)\right)= \\
&=&(-1)^{n+1} \left( \sum_{j=1}^n j S_{n,j}
(t+r-1)^j\frac{\partial^j}{\partial t^j}E_{r,0}(t) +
\sum_{j=2}^{n+1} (j-1) S_{n,j-1} (t+r-1)^j\frac{\partial^j}{\partial
t^j}E_{r,0}(t)\right)=\\
&=&(-1)^{n+1} \sum_{j=1}^{n+1} S_{n+1,j}
(t+r-1)^j\frac{\partial^j}{\partial t^j}E_{r,0}(t)
\end{eqnarray*}
\end{tiny}
as needed.

\medskip

Using the initial condition of this recurrence, namely
$E_{r,0}(t)=\frac{r}{t-1}$, we obtain that:
$$E_{r,n}(t)=(-1)^n\sum_{j=1}^nS_{n,j}(t+r-1)^j\frac{(-1)^jj!r}{(t-1)^{j+1}},$$
which is equivalent to:
$$E_{r,n}(t)=(-1)^{n-1}r\sum_{j=1}^{n}j!S_{n,j}\frac{(t+r-1)^j}{(1-t)^{j+1}},$$
which completes the proof.
\end{proof}

Now we ready to give an explicit formula for the polynomial
$D_{r,n}(t)$.

\begin{cor}
For all $n\geq1$,
$$D_{r,n}(t)=r\sum_{j=1}^nj!S_{n,j}(t+r-1)^{j-1}(1-t)^{n-j}.$$
\end{cor}
\begin{proof}
From the definitions, we have that
$$D_{r,n}(t)=P_{r,n}(t)E_{r,n}(t)=\frac{(t-1)^{n+1}}{t+r-1}E_{r,n}(t)=\frac{(-1)^{n+1}(1-t)^{n+1}}{t+r-1}E_{r,n}(t).$$
Now, by Theorem \ref{thmdd5}, we get the desired result.
\end{proof}

\begin{rem}
The solution of this type of a recurrence relation is based on
\cite{MS}.
\end{rem}

 Finding the coefficient of $t^k$ in the polynomial
$D_{r,n}(t)$ we obtain an explicit formula for the element
$d(r,n,k)$, as follows.

\bthm The number of permutations $\pi$ in $\grn$ with satisfy ${\rm
exc}_A(\pi)=k$ is given by
$$d(r,n,k)=r\sum_{j=1}^n\sum_{i=0}^{j-1}(-1)^{k+j-1-i}r^ij!S_{n,j}\binom{j-1}{i}\binom{n-1-i}{k}.$$
\ethm

\section{Log-concavity of the parameter ${\rm exc}_A$}\label{concave}

In this section, we show that the parameter ${\rm exc}_A$ is
log-concave. We prove that ${\rm exc}_A$ is log-concave in $B_n$.
The corresponding proof for $\grn$ is similar.

Since log-concavity implies unimodality for positive sequences, the
parameter ${\rm exc}_A$ is unimodal too.

\bthm The parameter ${\rm exc}_A$ on $B_n$ is log-concave. \ethm

\begin{proof}
Recall the following definitions from Section \ref{proof_Bn}:
$$C_i(n,k)=\{ \pi \in B_n \mid \exc_A(\pi)=k, \negg(\pi)=i \},$$
$$c_i(n,k)=\# C_i(n,k).$$

Denote: $$D_{n,k}=\sum_{i=0}^{n}{c_i(n,k)}.$$

By the recursion given in Proposition \ref{Cink}:
\begin{eqnarray*}
c_i(n,k)&=&(n-k)c_i(n-1,k-1)+ (k+1)c_i(n-1,k)+ \\
& &+(n-k)c_{i-1}(n-1,k)+(k+1)c_{i-1}(n-1,k+1),
\end{eqnarray*}
we have:
\begin{small}
\begin{eqnarray*}
D_{n,k}&=&(n-k)D_{n-1,k-1}+(k+1)D_{n-1,k}+(n-k)D_{n-1,k}+(k+1)D_{n-1,k+1}= \\
& =&(n-k)D_{n-1,k-1}+(n+1)D_{n-1,k}+(k+1)D_{n-1,k+1}
\end{eqnarray*}
\end{small}

We prove the log-concavity by induction. For $n=3$, the claim can be easily verified. Now we assume it for $n-1$, and we have to show that:
$$D_{n,k}^2 \geq D_{n,k-1} D_{n,k+1}.$$

Along the following computation, we abbreviate $D_{n-1,j}$ to $D_j$.

We compute:

\begin{tiny}
\begin{eqnarray*}
D_{k}^2 - D_{k-1} D_{k+1}& =&  \left[ (n-k)D_{k-1}+(n+1)D_k+(k+1)D_{k+1} \right]^2 - \\
& & - \left[ (n-k+1)D_{k-2}+(n+1)D_{k-1}+k D_k \right] \cdot \\
& & \qquad \cdot \left[ (n-k-1)D_k+(n+1)D_{k+1}+(k+2)D_{k+2} \right]= \\
& & \\
&=& \left[ (n-k)^2 D^2_{k-1} + (n+1)^2 D_k^2 + (k+1)^2 D^2_{k+1} +  \right. \\
& & \qquad \left. + 2(n-k)(n+1) D_{k-1} D_k +2(n-k)(k+1)D_{k-1} D_{k+1} + \right. \\
& & \qquad \left. + 2(n+1)(k+1) D_k D_{k+1} \right] - \\
& & - \left[ (n-k+1)(n-k-1)D_{k-2}D_k +(n-k+1)(n+1)D_{k-2}D_{k+1}+  \right.  \\
& & \qquad + \left. (n-k+1)(k+2)D_{k-2}D_{k+2} + (n+1)(n-k-1)D_{k-1}D_k +\right. \\
& & \qquad + \left. (n+1)^2 D_{k-1}D_{k+1} + (n+1)(k+2)D_{k-1}D_{k+2} +\right. \\
& & \qquad + \left. k(n-k-1) D_k^2 + k(n+1)D_kD_{k+1} +k(k+2)D_k D_{k+2} \right]= \\
& & \\
&=& \left[ (n-k)^2 D^2_{k-1} - (n-k+1)(n-k-1)D_{k-2}D_k \right]+\\
& & + \left[ (k+1)^2 D^2_{k+1} - k(k+2)D_k D_{k+2} \right]+\\
& & + \left[ 2(n-k)(n+1) D_{k-1} D_k - (n-k+1)(n+1)D_{k-2}D_{k+1} - \right. \\
& & \qquad - \left. (n+1)(n-k-1)D_{k-1}D_k \right]+\\
& & + \left[ 2(n+1)(k+1) D_k D_{k+1} - (n+1)(k+2)D_{k-1}D_{k+2}  - \right. \\
& & \qquad - \left. k(n+1)D_kD_{k+1} \right]+ \\
& & +\left[ (n+1)^2 D_k^2  + 2(n-k)(k+1)D_{k-1} D_{k+1} -  \right. \\
& & \qquad - \left. (n-k+1)(k+2)D_{k-2}D_{k+2} - (n+1)^2 D_{k-1}D_{k+1} - \right. \\
& & \qquad - \left. k(n-k-1) D_k^2 \right]\\
\end{eqnarray*}
\end{tiny}

We treat each one of the five brackets separately.

By the induction hypothesis, the first bracket is greater (or equal) than $D^2_{k-1}$.
Similarly, the second bracket is greater (or equal) than $D^2_{k+1}$.

Since $D_{k-2}D_{k+1} \leq D_{k-1}D_k$ by the log-concavity assumption, we have:

\begin{tiny}
$$2(n-k)(n+1) D_{k-1} D_k - (n-k+1)(n+1)D_{k-2}D_{k+1} - (n+1)(n-k-1)D_{k-1}D_k \geq $$
$$\geq 2(n-k)(n+1) D_{k-1} D_k - (n-k+1)(n+1)D_{k-1}D_k - (n+1)(n-k-1)D_{k-1}D_k=0 \cdot D_{k-1}D_k =0,$$
\end{tiny}

\noindent
and hence the third bracket is non-negative. Similarly, the
fourth bracket is non-negative.

From the first two brackets, we have two positive elements:
$D^2_{k-1}$ and $D^2_{k+1}$. Their sum can be written as
$(D_{k-1}-D_{k+1})^2+2D_{k-1}D_{k+1}$. Adding $2D_{k-1}D_{k+1}$ to
the fifth bracket, we have:

\begin{tiny}
$$(n+1)^2 D_k^2  + 2((n-k)(k+1)+1)D_{k-1} D_{k+1} - (n-k+1)(k+2)D_{k-2}D_{k+2}-$$
$$ - (n+1)^2 D_{k-1}D_{k+1} - k(n-k-1) D_k^2$$
\end{tiny}

Since $(n+1)^2=((n-k)+(k+1))^2=(n-k)^2+(k+1)^2+2(n-k)(k-1)$, we can
simplify the bracket into:
$$(n+1)^2 D_k^2 - (n-k+1)(k+2)D_{k-2}D_{k+2}-$$
$$ - ((k+1)^2-1) D_{k-1}D_{k+1} - ((n-k)^2-1) D_{k-1}D_{k+1}- k(n-k-1) D_k^2$$
By the log-concavity assumption, this sum is greater (or equal) than:
\begin{tiny}
$$(n+1)^2 D_k^2 - (n-k+1)(k+2)D_k^2- ((k+1)^2-1) D^2_k - ((n-k)^2-1) D_k^2- k(n-k-1) D_k^2$$
\end{tiny}
which is equal to $0$, and hence we show that the sum of all five brackets is non-negative and hence we are done.
\end{proof}

\section{Symmetry of the excedance statistic}

In this section, we present the following symmetry property of the
generating function of the excedance statistics on $G_{r,n}$:

\bthm The polynomial $$\sumlim_{\pi \in
\grn}{{q^{\exc(\pi)}}}=\sumlim_{i=0}^{rn-1}{a_iq^i}$$ satisfies:
$a_i=a_{rn-1-i}$ for $i \in \{1,\cdots,\lfloor \frac{rn-1}{2}
\rfloor\}$.

\ethm

Define a bijection of $\grn$: $\pi \mapsto \pi'$ in the following
way: \\ For $1 \leq i \leq n-1$:
\begin{itemize}
\item If $\pi(i)=j^{[\beta]}$ then $\pi'(n-i)=(n+1-j)^{[r-\beta]}$.
\item If $\pi(n)=j^{[\beta]}$ then $\pi'(n)=(n+1-j)^{[r-1-\beta]}$.
\end{itemize}

Instead of burdening the reader with the subtle though standard
proof, we choose to give an example.

\begin{exa}

$$\pi=\left(\begin{array}{cccc|cccc|cccc}
\bar{\bar{1}} & \bar{\bar{2}} & \bar{\bar{3}} & \bar{\bar{4}} & \bar{1} & \bar{2} & \bar{3} & \bar{4} & 1 & 2 & 3 & 4 \\
2 & \bar{1} & 4 & \bar{\bar{3}} & \bar{\bar{2}} & 1 & \bar{\bar{4}} &  \bar{3} & \bar{2} & \bar{\bar{1}} & \bar{4} &3 \\
\end{array}\right)
$$

$$\pi'=\left(\begin{array}{cccc|cccc|cccc}
\bar{\bar{1}} & \bar{\bar{2}} & \bar{\bar{3}} & \bar{\bar{4}} & \bar{1} & \bar{2} & \bar{3} & \bar{4} & 1 & 2 & 3 & 4 \\
\bar{1} & 4 & \bar{3} & 2 & 1 & \bar{\bar{4}} &  3 & 2 & \bar{\bar{1}} & \bar{4} & \bar{\bar{3}} & \bar{\bar{2}} \\
\end{array}\right).
$$

The important point, one can easily see here, is that in any place
there is an excedance at $\pi$, there is no excedance at $\pi'$, and
vice versa.
\end{exa}

\end{document}